\documentclass[runningheads]{llncs}
\usepackage{enumerate}
\usepackage{ifthen,xspace}
\usepackage{graphicx,color,MnSymbol}
\usepackage{amsmath,amsfonts,amstext,booktabs}
\makeatletter\@ifclassloaded{llncs}{}{\usepackage{amsthm}}\makeatother

\newtheorem{prop}{Proposition}

\newtheorem{lem}[prop]{Lemma}
\newtheorem{defn}[prop]{Definition}

\newtheorem{thm}[prop]{Theorem}

\newcommand{\Ord}{\mathrm{Ord}}
\newcommand{\cOrd}{\Sigma^\mathrm{O}}

\newcommand{\nothing}{{\boldsymbol{\circ}}}

\newcommand{\GBC}{\mathsf{NBGC}}
\newcommand{\GB}{\mathsf{NBG}}
\newcommand{\AC}{\mathsf{AC}}
\newcommand{\ZFC}{\mathsf{ZFC}}

\newcommand{\leftsym}{\filledmedtriangleleft}
\newcommand{\rightsym}{\filledmedtriangleright}
\newcommand{\staysym}{\filledmedtriangledown}

\newcommand{\VeqL}{{\mathrm{V}{=}\mathrm{L}}}

\newcommand{\zero}{\mathbf{0}}
\newcommand{\one}{\mathbf{1}}

\newcommand{\MT}{\mathrm{MT}}

\newcommand{\bfm}{\mathbf{m}}

\newcommand{\converges}{{\downarrow}}
\newcommand{\diverges}{{\uparrow}}
\newcommand{\dom}{\mathrm{dom}}
\newcommand{\ran}{\mathrm{ran}}
\newcommand{\pto}{\dashrightarrow}

\newcommand{\infwords}{\Sigma^{{<}\mathrm{Ord}}}
\newcommand{\rinfwords}{\Sigma^{(\mathrm{Ord})}}
\newcommand{\snap}{\mathrm{Snap}}

\newcommand{\FS}{\Sigma^\mathrm{FS}}

\newenvironment{proofsketch}{\noindent\emph{Proof sketch.}}{\hfill\textsc{q.e.d.}\\
\smallskip}

\newcommand{\wech}[1]{}

\begin{document}

\allowdisplaybreaks

\authorrunning{L.\ Galeotti, E.\ S.\ Lewis, B.\ L\"owe}
\titlerunning{Symmetry for transfinite computability}

\author{Lorenzo Galeotti$^1$, Ethan S.\ Lewis$^{2}$, Benedikt L\"owe$^{3,4,5}$}

\title{Symmetry for transfinite computability}

\institute{
Amsterdam University College, Postbus 94160,
1090 GD Amsterdam, The Netherlands; \texttt{l.galeotti@uva.nl}
\and
1179 W 1200 S,
Springville, UT 84663,
United States of America;
\texttt{etonlels@gmail.com}
\and
Institute for Logic, Language
and Computation, Universiteit van Amsterdam, Postbus 94242, 1090\,GE~Amsterdam, The Netherlands; \texttt{b.loewe@uva.nl} \and
Fachbereich Mathematik, Universit\"at Hamburg, Bundesstrasse
55, 20146 Hamburg, Germany
\and
Churchill College, Lucy Cavendish College, St.\ Edmund's College, \& Department of Pure Mathematics and Mathematical Statistics,
University of Cambridge,
Storey's Way,
Cambridge, CB3 0DS, England
}

\maketitle

\begin{abstract}
Finite Turing computation has a fundamental symmetry between inputs, outputs, programs, time, and storage space.
Standard models of transfinite computational break this symmetry; we consider ways to recover it and study the resulting model of computation.
This model exhibits the same symmetry
as finite Turing computation in universes constructible from a set of ordinals, but that statement is independent of 
von Neumann-G\"odel-Bernays class theory.
\end{abstract}

\setcounter{footnote}{0}

\section{Introduction}
\label{sec:intro}

A fundamental feature of the theory of computation is that the constituents of computability, viz.\
in-/output, programs, time, and storage space can be considered to be the same type of object: natural numbers (if necessary, via coding). 
A Turing machine receives a finite string of symbols as input, has a finite string of symbols as its program, and produces a finite string of symbols as output. Moreover, both its tape and its time flow are indexed by natural numbers. Therefore, since finite strings of symbols can be coded as a natural number, all these objects are of the same type.

We shall refer to this feature as \emph{symmetry}. Various aspects of symmetry permeate the
general theory of computation: the symmetry between inputs and programs is the reason for the
\emph{software principle} (the 
existence of universal machines) and the $s$-$m$-$n$ Theorem; the symmetry between programs and time underlies the zigzag method that allows us to parallelise infinitely
many computations into one by identifying the cartesian product of the program space and time with $\mathbb{N}\times\mathbb{N}$ and using Cantor's zigzag function.

The oldest model of transfinite computation are the \emph{Hamkins-Kidder machines} or \emph{Infinite Time Turing Machines} (ITTM), defined in \cite{HamLew00}.
These machines have a \emph{storage space} of order type $\omega$, but allow computation to be of arbitrary transfinite ordinal length, thereby breaking the symmetry between time and space. This asymmetry makes their complexity theory
vastly different from ordinary complexity theory, as discussed in
\cite{schindler,hamwel,deohamsch,loewe2006,winterthesis,winter2009}.

In \cite{koepke2005,koepke2009}, Koepke symmetrised Hamkins-Kidder machines and defined what is now known as
\emph{Koepke machines} or \emph{Ordinal Turing Machines}: Koepke machines have a class-sized tape indexed by ordinals and run through ordinal time, thereby re-establishing the symmetry between \emph{time} and \emph{storage space}.\footnote{Carl argues in  \cite[Chapter 9]{carl} that Koepke machines are the natural infinitary analogue for finitary computation and complexity theory and this was explored in detail in \cite{CLR}.}
However, Koepke machines do not have the full symmetry
that we find in finite Turing computation: while \emph{time} and \emph{storage space} are represented by arbitrary ordinals, programs are still finite objects.

In this paper,
we shall provide a general framework for models of computation and computability that allows us to phrase the quest for symmetry in abstract terms; this is done in \S\,\ref{sec:framework}.
In this framework, we shall define the relevant models of computability, i.e.,
ordinary Turing computability, Hamkins-Kidder computability, Koepke computability, and our new
notion called \emph{symmetric computability} in \S\,\ref{sec:models}. We study basic properties of
symmetric computability in \S\,\ref{sec:treble}, and finally show that the full symmetry of symmetric computability cannot be proved in von Neumann-G\"odel-Bernays class theory ($\GBC$) in \S\,\ref{sec:ce}: symmetry holds if and only if the universe is constructible from a set of ordinals.

This paper contains results from the second author's Master's thesis \cite{lewis}
written under the supervision of the first and the third author. These results are
cited in 
\cite{carl2022,Galeotti2019CiE,Carl2020} and
\cite[Exercise 3.9.7]{carl}.

\section{Class theories}\label{sec:classtheory}

In this paper, we work in \emph{von Neumann-G\"odel-Bernays class theory}.\footnote{For more details, cf., e.g., \cite[Chapter 4]{mendelson}.}
The language is the usual language of set theory $\mathcal{L}_\in$ with a single binary relation symbol~$\in$. We define a unary predicate $\mathrm{set}(x) := \exists y (x\in y)$. Using this predicate, we can define the two \emph{set quantifiers}
$\exists^\mathrm{set}x \varphi := \exists x (\mathrm{set}(x)\wedge \varphi)$ and
$\forall^\mathrm{set}x \varphi := \forall x (\mathrm{set}(x)\to \varphi)$. A formula is called \emph{set theoretic} if all of its quantifiers are set quantifiers.
In this context we 
denote by $\AC$ the \emph{axiom of choice} for sets, i.e., the statement that  ``Every set $x$ has a choice function'' and contrast it with the axiom of \emph{Global Choice} which is the statement ``There is a global choice class function''. We write $\GB$ for von Neumann-G\"odel-Bernays class theory without the axiom of Global Choice
\cite[p.\ 70: Axioms A--D]{jech}
and $\GBC$ for the theory obtained from $\GB$ adding the axiom of Global Choice
\cite[p.\ 70: Axioms A--E]{jech}.
It is a well-known result due to Easton that 
if 
$\GB$ is consistent, then
$\GB+\AC$ does not prove the axiom of Global Choice
(cf., e.g., {\cite[Theorem 3.1]{felgner}}).

We can transform a formula $\varphi$ in the language $\mathcal{L}_\in$ into a set theoretic formula $\varphi^\mathrm{set}$ by recursively replacing all quantifiers with the corresponding set quantifiers.
This allows us to formulate the famous \emph{conservativity theorem for $\GBC$} (cf., e.g., \cite[Corollary 4.1 \& Theorem 4.2]{felgner}):\footnote{We refer the reader to \cite[p.\ 242]{felgner} and \cite[p.\ 381]{ferreiros} for more information on the history of this theorem.}

\begin{thm}\label{thm:cons_ext}
If $\varphi$ is a sentence in $\mathcal{L}_\in$, then $\ZFC\vdash\varphi$ if and only if
$\GBC\vdash \varphi^\mathrm{set}$ if and only if
$\GB+\AC\vdash \varphi^\mathrm{set}$.
\end{thm}



We define the 
\emph{axiom of constructibility from a set of ordinals} as the statement ``there is a set of ordinals $x$ such that
$\mathrm{V}{=}\mathrm{L}[x]$''.
This 
is a set theoretic sentence and 
implies the axiom of Global Choice; thus Theorem 
\ref{thm:cons_ext} implies the following result.

\begin{thm}\label{thm:GBC_V=L} If $\GBC$ is consistent, then
$\GBC$ does not prove nor disprove the axiom
of constructibility from a set of ordinals.
\end{thm}

\section{The general framework of Turing computation and computability}
\label{sec:framework}

We shall frame our discussion of symmetry in a general context that makes the relavent models
of computability special cases of a general framework. Our general framework will work on the class $\Ord$
of all ordinals and refer to the class $\mathrm{V}$ of all sets as potential programs for these machines.
All models of computation in this paper will be variants of Turing machines: they have a single 
class-length tape indexed by ordinals,\footnote{For most models of computability, the number of tapes does not matter; however, in the case of Hamkins-Kidder machines,
1-tape machines and 3-tape machines differ (cf.\ \cite{HamSea}). Since we do not discuss Hamkins-Kidder machines in detail, this is immaterial for our context. \label{fn:3tape}}
a read/write head that moves on the tape according to a program.
We fix a finite alphabet $\Sigma$ with at least two elements $\zero$ and $\one$ for the remainder of the paper.

\medskip

\noindent\textbf{Turing hardware \& computations.} At the highest level of abstraction, we deal with the \emph{Turing hardware}: the \emph{tape} and the \emph{head}, including the description of how they work. We assume that the tape is always indexed by ordinals, split up into discrete \emph{cells} in which a symbol from $\Sigma$ can be written; also, we assume
that time is considered as discrete points in time, indexed by ordinals, and that at each point in time, the head is located at one of the cells; finally, we assume that we have discrete \emph{states}, indexed by ordinals. 

For classes $X$ and $Y$, we write
$f\colon X\dashrightarrow Y$ for ``$f$ is a class function
with $\dom(f)\subseteq X$ and $\ran(f)\subseteq Y$'' 
and 
$f(x)\converges$ if and only if $x\in\dom(f)$ and $f(x)\diverges$ otherwise. We represent the tape content by arbitrary
partial class functions from $\Ord$ to $\Sigma$; we write $\rinfwords$ for the class of these objects.
We shall consider a number of relevant subclasses of this class: $\infwords := \{x\in\rinfwords\,;\,\mathrm{dom}(x)\in\Ord\}$,
 $\FS := \{x\in\rinfwords\,;\,\mathrm{dom}(x)\mbox{ is finite}\}$,
$\Sigma^\omega := \{x\in\rinfwords\,;\,\mathrm{dom}(x)=\omega\}$,
$\Sigma^* = \Sigma^{<\omega} := \{x\in\rinfwords\,;\,\mathrm{dom}(x)\in\omega\}$, and
$\cOrd  := \{x\in\rinfwords\,;\,|\mathrm{dom}(x)| = 1\mbox{ and }\mathrm{ran}(x) = \{\zero\}\}$.

The classes $\cOrd$ and $\FS$ are our representations of the classes
$\Ord$ and $\Ord^{<\omega}$, respectively. The classes $\Ord$
and $\cOrd$ have a canonical bijection; the classes $\FS$ and $\Ord^{<\omega}$
can be identified via the G\"odel pairing function.\footnote{The Gödel pairing function is an absolutely definable class bijection between
$\mathrm{Ord}$ and $\mathrm{Ord}^2$; cf.\ \cite[pp.\ 30--31]{jech}.}
Furthermore, the Gödel pairing function yields a definable bijection between
$\Ord$ and $\Ord^{<\omega}$ and 
a bijection $(w,v)\mapsto w*v\colon
\infwords\times\infwords\to\infwords$.

A \emph{snapshot} of the machine consists of the tape content, a state, and a position of the head, i.e.,
a tuple from $\snap:= \rinfwords\times\Ord\times\Ord$.

The behaviour of the head is governed by the \emph{transition rule},
a class function that describes what the head will do
given its past behaviour and a program $p$. 
For now, we still allow all sets to be programs
(we consider this specification to be part of the \emph{software}), so a transition function is a
class function $T\colon\snap^{{<}\Ord}\times\mathrm{V}\to\snap$.
Once a transition function $T$ is fixed, 
given a program $p$ and a snapshot 
$s = (x,\alpha,\beta)$, we 
define an ordinal-length sequence of snapshots by recursion:
$C_{p,s}(0):= s$, and
$C_{p,s}(\gamma):= T(C_{p,s}{\upharpoonright}\gamma,p)$ for $\gamma>0$.
We shall call this \emph{the computation of program $p$ with initial snapshot $s$}.

In this paper, we shall only consider two different transition functions,
the \emph{finite transition function} $T_\mathrm{f}$ which is used by ordinary Turing machines and
Hamkins-Kidder machines and the \emph{transfinite transition function} $T_\mathrm{t}$ which is used by Koepke machines (for definitions, cf.\ \S\,\ref{sec:models}).

\medskip

\noindent\textbf{Turing software.} A \emph{model of computation} consists of hardware (i.e., a transition function $T$) and
a class of programs $P$ that can be used for computing. Specifying the class of programs
identifies which of the computations are computations according to
a program in $P$.

In this paper, we shall only consider two classes of programs, the class of \emph{finite programs}
$P_\mathrm{f}$ and the class of \emph{transfinite programs} $P_\mathrm{t}$ (for definitions, cf.\ \S\,\ref{sec:models}).

\medskip

\noindent\textbf{Computability.} A model of computation determines a class of computations, but does not yet tell us what they do. To illustrate this, consider
the ordinary notion of Turing computation: 
for each program and snapshot, we get an infinite sequence of snapshots, but there are many ways to interpret these infinite sequences. Following Turing's original seminal definition \cite[\S\,2]{turing}, we designate 
\emph{start} and \emph{halt states}, give a definition of \emph{halting computations} and then interpret the computation as producing a partial function (for definitions, cf.\ \S\,\ref{sec:models}).

Abstractly, we say that an \emph{interpretation}
consists of 
a partial class function $I$
that assigns to a transition function $T$ and each program $p\in P$ a partial
function $I(T,p)\colon \rinfwords\dashrightarrow \rinfwords$
and a class $D\subseteq\rinfwords$ called
the \emph{domain} of the interpretation.
A \emph{model of computability} is a model of computation (i.e., a transition function $T$ and
a class of programs $P$) together with an interpretation. We say that $f\colon D\dashrightarrow D$ is
\emph{computable} according to this model of computability if there is a $p\in P$ such that
$f = I(T,p){\upharpoonright}D$.

Note that for a given model of computation and a fixed interpretation function, there is some freedom to choose $D$.
E.g., usually, for ordinary Turing computations with the usual textbook interpretation,
we let $D = \Sigma^*$ and
thus, 
computability is a property of partial functions $f\colon \Sigma^*\dashrightarrow \Sigma^*$.
However, we could consider $D = \Sigma^\omega$, i.e., letting the Turing machine operate on arbitrary
tape contents of length $\omega$, obtaining a different model of computability.\footnote{This is a curious model of computability that exhibits a discrepancy between 1-tape and 3-tape machines
(cf.\ Footnote \ref{fn:3tape}): 
since only finitely many cells are changed in halting
computations, 
for 1-tape machines the identity function is  computable and constant functions are not; in contrast, for 3-tape machines 
constant functions with value $w\in\Sigma^*$ are computable, but other constant functions or the identity function are not. \label{fn:weird}}
On the other hand, if you fix the model of computation and the type of interpretation function, $D$ cannot be chosen entirely freely: the class $D$ needs to be closed under the operations $I(T,p)$ for $p\in P$.
E.g., if our model of computation is Koepke machines with the usual interpretation, we cannot choose $D = \Sigma^*$ or even $D=\Sigma^\omega$ since there are
programs with which a Koepke machine would produce an output that is not in $D$ anymore.

In this paper, we shall consider two types of interpretation function,
the \emph{finite interpretation} $I_\mathrm{f}$ and
the \emph{transfinite interpretation} $I_\mathrm{t}$
(for definitions, cf.\ \S\,\ref{sec:models}).

\section{Concrete models of computation}\label{sec:models}

\noindent\textbf{Programs.} 
We fix three \emph{motion tokens}
$\MT := \{\leftsym,\staysym,\rightsym\}$ 
that represent the instructions for the head movements (``move left'', ``do nothing'', and ``move right'') and use $\bfm$ as variable for motion tokens.
Among the states (indexed by ordinals), we single out three particular states: the \emph{start state} indexed by $0$,
the \emph{halt state} indexed by $1$, and the \emph{limit state} indexed by $2$.
We write $\Sigma_\nothing := \Sigma\cup\{\nothing\}$
where $\nothing$ is a special symbol representing an empty cell.
All of our programs will be
partial functions $p\colon \Ord\times\Sigma_\nothing\pto \Ord\times\Sigma_\nothing\times\MT$.
We call a program \emph{finite} if its domain is finite and \emph{transfinite} if its domain is a set.
The classes of finite and transfinite programs are denoted by $P_\mathrm{f}$ and $P_\mathrm{t}$,
respectively.

Via the canonical identification of the classes
$\Ord$, $\Ord\times\Sigma_\nothing$, and $\Ord\times\Sigma_\nothing\times\MT$,
we can encode programs as elements of 
$\rinfwords$. 
Under our encoding, we identify 
$P_\mathrm{f}$ with the class $\FS$ and $P_\mathrm{t}$ with the class $\infwords$.

\medskip

\noindent\textbf{Transition functions.} Given a program $p$,
we shall now define the transition functions $T_\mathrm{f}$ (``finite transition function'')
and $T_\mathrm{t}$ (``transfinite transition function''). They are identical
on sequences of successor length and coincide there with the ordinary
transition function defined by Turing for his machines; they differ for
sequences of limit length.

If $\vec s = (s_\xi\,;\,\xi<\gamma+1)$ is a sequence of snapshots of successor
length, the transition function will only depend on $s_\gamma = (x,\alpha,\beta)$,
the final snapshot in the list. Thus, $x\in\rinfwords$ is the tape content at time $\gamma$,
$\alpha$ is the state at time $\gamma$, and $\beta$ is the location of the head at time $\gamma$.
If $p(\alpha,x(\beta))$ is undefined, we let $T(\vec s) := s_\gamma$; otherwise, let
$p(\alpha,x(\beta)) = (\delta,\sigma,\bfm)$. Then $T(\vec s) = (y,\alpha^+,\beta^+)$
where $\alpha^+  := \delta$, 
\smallskip

\begin{minipage}{.32\textwidth}
$y(\eta):= \begin{cases}
x(\eta) & \mbox{ if $\eta \neq \beta$,}\\
\sigma & \mbox{ if $\eta = \beta$,}
\end{cases}$
\end{minipage} and %
\begin{minipage}{0.5\textwidth}
$\beta^+ :=\begin{cases}
\beta-1 & \mbox{ if $\bfm = \leftsym$ and $\beta$ is a successor,}\\
0 & \mbox{ if $\bfm = \leftsym$ and $\beta$ is a limit,}\\
\beta+1 & \mbox{ if $\bfm = \rightsym$,}\\
\beta & \mbox{ if $\bfm = \staysym$.}
\end{cases}$        
\end{minipage}

\smallskip

If $\vec s = (s_\xi\,;\,\xi<\lambda)$
with $s_\xi = (x_\xi,\alpha_\xi,\beta_\xi)$
is a sequence of snapshots of
limit length $\lambda$, the two transition functions agree in their definition of
the tape content, but disagree in their treatment of the head position and
state. Let us write
$T_\mathrm{f}(\vec s) = (y,\alpha_\mathrm{f},\beta_\mathrm{f})$
and
$T_\mathrm{t}(\vec s) = (y,\alpha_\mathrm{t},\beta_\mathrm{t})$.
For the tape content, we assume that we have a total ordering on $\Sigma$
and define
$y(\eta) := \mathrm{liminf}\{x_\xi(\eta)\,;\,\xi<\lambda\}$.

The finite transition function $T_\mathrm{f}$ moves the head to cell $0$,
moves to the limit state (indexed by $2$), i.e.,
$\alpha_\mathrm{f} := 2$ and $\beta_\mathrm{f} := 0$. Note that in any computation using the finite transition function, the head will never
reach a cell indexed by an infinite ordinal.

The transfinite transition function $T_\mathrm{t}$ moves both the head and the
cell to the \emph{inferior limit} of the ordinals occurring in the sequence, i.e.,
$\alpha_\mathrm{t} :=
\mathrm{liminf}\{\alpha_\xi\,;\,\xi<\lambda\}$
and
$\beta_\mathrm{t} :=
\mathrm{liminf}\{\beta_\xi\,;\,\xi<\lambda ~\land~ \alpha_\xi= \alpha_{\mathrm{t}} \}$.

\medskip

\noindent\textbf{Interpretations.}
We define our two interpretation functions uniformly for arbitrary tape contents $x\in\rinfwords$.
Both interpretations take a tape content $x$ and a program $p$, define the initial
snapshot $s := (x,0,0)$, and produce the computation $C_{p,s}$ of program $p$ with initial snapshot $s$.

The \emph{finite interpretation} $I_\mathrm{f}$ considers a computation as \emph{halting} if there is a natural number $n$ such that the state of $C_{p,s}(n)$ is $1$ (i.e., the halting state); the \emph{transfinite interpretation} $I_\mathrm{t}$
considers a computation as \emph{halting} if there is an ordinal $\alpha$ such that the state of
$C_{p,s}(\alpha)$ is $1$. If it exists, the smallest such number is called the \emph{halting time}
of the computation.
This implicitly defines the time considered by these models of computability:
in general, we say that the \emph{time relevant for a model of computability} is the 
supremum of its halting times. This is at most $\omega$ for models with the finite interpretation
and at most $\Ord$ for models with the infinite interpretation. We let $\Omega\subseteq D$ be a subclass that is identified with the time relevant of the model, e.g., $\Omega = \cOrd$ if the relevant time is $\Ord$.

If a computation is halting, we say that the tape content at its halting time is the \emph{output}
of the computation.
Finally, for $I = I_\mathrm{f}$ or $I = I_\mathrm{t}$ and the appropriate notion of
halting, we let $I(T,p)(x) := y
$ if $C_{p,s}$ is halting and $y$ is its output, and $I(T,p)(x)\diverges$  otherwise.

\medskip

\noindent\textbf{Models of computability.}
Using our finite specifications $T_\mathrm{f}$, $P_\mathrm{f}$, and $I_\mathrm{f}$
and our transfinite specifications $T_\mathrm{t}$, $P_\mathrm{t}$, and $I_\mathrm{t}$,
we can now recover the known models of computability (and a new one) as special cases.

First observe that if the transition function is finite, then any tape content beyond the cells indexed by natural numbers will be immaterial for the computation since the head never moves to these cells. So, the
relevant input domain has to be $\Sigma^*$ or
$\Sigma^\omega$. Moreover, if the interpretation function is finite, then all computations that go on to
$\omega$ or beyond will be disregarded in the interpretation, so we can assume, without loss of generality,
that the transition function is finite as well. This leads to the models of computability listed in
Table \ref{table:1}.

\begin{table}[h]
\begin{center}{\footnotesize \begin{tabular}{cccccl}
\toprule
& Transition & Programs & Interpretation & In-/Output & \\
\midrule
(a) & Finite & Finite & Finite & $\Sigma^*$ & ordinary Turing machines\\
(b) & Finite & Finite & Transfinite & $\Sigma^\omega$ & Hamkins-Kidder machines\\
(c) & Transfinite & Finite & Transfinite & $\infwords$ & Koepke machines\\
(d) & Transfinite & Transfinite & Transfinite & $\infwords$ & symmetric machines; cf.\ \S\,\ref{sec:treble}\\
\bottomrule
\end{tabular}}\end{center}
\caption{The considered models of computability\label{table:1}}
\end{table}

We briefly discuss the choice of in-/output for the described models:

\smallskip

\noindent \emph{Table \ref{table:1} (a).} Since the transition function
is finite, only $\Sigma^*$ and $\Sigma^\omega$ make sense as choice of in-/output.
However, since the interpretation is
also finite, no halting computation will ever be able to read an entire infinite
tape, so $\Sigma^*$ is the natural choice for in-/output. Choosing
$\Sigma^\omega$ leads to the model of computability
discussed in Footnote \ref{fn:weird}.
The time relevant for this model of $\omega$.

\smallskip

\noindent \emph{Table \ref{table:1} (b).} Similarly in this case, the finite transition
function means that we can only choose $\Sigma^*$ or $\Sigma^\omega$ as input;
however, $\Sigma^*$ is not closed under the operation of the interpretation
(a Hamkins-Kidder machine can fill the entire tape and then halt), so $\Sigma^\omega$ is the only remaining natural choice.
The time relevant for Hamkins-Kidder
machines
has been investigated in \cite{HamLew00,welch}.

\smallskip

\noindent \emph{Table \ref{table:1} (c).} In analogy to the argument given for
line (a), a halting Koepke machine will only consider a set of cells
on the tape. Thus, the natural choice of in-/output is $\infwords$. Similarly to
(a), it makes sense to consider $D = \rinfwords$
in which case the discussion of Footnote~\ref{fn:weird} applies.
The time relevant for Koepke computability is the class $\Ord$ of all ordinals.
If $p\in\infwords$, we say that a partial function $f\colon\infwords\dashrightarrow\infwords$ is \emph{Koepke computable with parameter $p$} if
the partial function $p*w\mapsto f(w)$  is Koepke computable. We shall prove in Proposition \ref{prop:symKoepke} that this notion is equivalent to the new notion introduced in line (d).

\smallskip

\noindent \emph{Table \ref{table:1} (d).}
The notion of computability introduced in line (d), called \emph{symmetric computability},
corrects the lack of symmetry
for transfinite computability.
In this model, time, space, and programs are all transfinite. The time relevant for symmetric computability is the class $\Ord$ of all ordinals.

\section{Symmetric machines}
\label{sec:treble}

In Table \ref{table:1}, we defined the \emph{model of symmetric computability} to be given by
the transfinite transition function, transfinite programs, and the transfinite
interpretation, using $\infwords$ as input and output. We call the corresponding
model of computation \emph{symmetric machines}.
In contrast to Hamkins-Kidder machines (who have considerably more time than space)
and Koepke machines (whose programs are tiny compared to the time and space they have available),
symmetric machines have set-sized 
time, space, and programs. They are the model of computability that systematically replaces the word ``finite'' in ordinary Turing computation with ``set-sized''.

\begin{prop}\label{prop:setsarecomputable}
If $f:\infwords\pto\infwords$ is a set, then $f$ is symmetrically computable.
\end{prop}

\begin{proofsketch}
Clearly, if $x\in\infwords$, there is a 
transfinite program $p$ that produces $x$ upon empty input (just explicitly specify the values of $x$). 

Since $f$ is a set, find some $\xi$ and 
$\Sigma^{(\xi)} := \{w\in\rinfwords\,\;\,\dom(w)\subseteq\xi\}$ such that both
$\dom(f)\subseteq\Sigma^{(\xi)}$ and $\ran(f)\subseteq\Sigma^{(\xi)}$.
By $\mathsf{AC}$, let $g$ be a bijection between 
some ordinal $\lambda$ and
$\Sigma^{(\xi)}$. 
We define $x_f\colon\lambda\cdot\xi\cdot 2\dashrightarrow\Sigma$
by letting $x_f(\alpha,\beta,0) := g(\alpha)(\beta)$
and $x_f(\alpha,\beta,1) := f(g(\alpha))(\beta)$.

Now find a program that writes $x_f$ on the tape.
Upon input $w\in\Sigma^{(\xi)}$, we can now determine
$f(w)$ as follows: search through the $0$-components of $x_f$ until you find $w$; when you found $w$ at index $\alpha$, 
output 
$\{(\beta,x_f(\alpha,\beta,1))\,;\,\beta < \xi\}$.
\end{proofsketch}

\begin{prop}\label{prop:symKoepke}
A partial function $f\colon\infwords\dashrightarrow\infwords$ is symmetrically
computable if and only of there is a $p\in\infwords$ such that $f$ is Koepke computable in parameter $p$.
\end{prop}

\begin{proofsketch}
By Proposition \ref{prop:setsarecomputable}, any 
parameter $p$ is symmetrically computable, 
so if $f$ is Koepke computable in parameter $p$, it is symmetrically computable as follows:
upon input $w$, first compute $p$, then $w*p$, then $f(p*w)$.

For the other direction, observe that in terms of hardware and interpretation, Koepke machines are just symmetric machines. Therefore,
the universal Koepke machine is also a universal symmetric machine, i.e., there is a Koepke machine $u$ such that for all transfinite programs $p$ and all $w\in\infwords$, we have
$$I_\mathrm{t}(T_\mathrm{t},u)(p*w) = I_\mathrm{t}(T_\mathrm{t},p)(w).$$ Thus, 
if a partial function $f$ is symmetrically computed by a program $p$, it is
Koepke computable with parameter $p$.
\end{proofsketch}

As usual, we can define the \emph{halting problem}
by
$$K := \{v*w\in\infwords\,;\,I_\mathrm{t}(T_\mathrm{t},v)(w){\downarrow}\}$$
(where $v$ is interpreted as a transfinite program).
The usual proof shows that $K$ is not symmetrically computable.

We shall now have a closer look at the symmetry properties for symmetric computability and ask whether it is the precise analogue of the symmetry exhibited by ordinary Turing computability.
For ordinary Turing machines, time and space are indexed by 
natural numbers; however, programs and in-/output are not \emph{prima facie} natural numbers;
they are finite sequences of elements of a finite set, i.e., via some encoding 
elements of $\Sigma^*$. In this case, the symmetry is given by the fact that
there is a computable encoding function that identifies $\Sigma^*$ and $\omega$.
Via such an encoding, we can see all four different parameters of the model of computability
as the same type of object.

In the case of symmetric computation, 
the word ``finite'' is systematically replaced by ``set-sized'', so time and space are indexed by ordinals and programs and in-/outputs are (up to encoding) elements of $\infwords$. 
Alas, in general, we cannot identify $\infwords$ and $\Ord$:
The encodings between the classes $\Ord$, $\Ord^{<\omega}$, $\cOrd$, and $\FS$ 
mentioned in \S\,\ref{sec:framework}
can be performed by Koepke machines (see, e.g., \cite[Section 4]{koepke2005}), but the class
$\infwords$ is a very different type of object: among other things, it contains the entire Cantor space (functions $f\colon \omega\to\Sigma$), so any computable encoding of elements $\infwords$ as ordinals would yield a computable wellordering of the reals. As a consequence, the existence of such a class function cannot be proved without additional set theoretic assumptions.

\section{The symmetry condition}
\label{sec:ce}

We give definitions of the notions of \emph{semidecidability} and \emph{computable enumerability}
within our abstract framework.
For the model of ordinary Turing computability, these definitions coincide with the usual definitions.

\begin{defn} 
Suppose that a model of computability is given by a transition function $T$, a class of programs $P$, and
an interpretation $I$ with domain class $D$.
Let $A\subseteq D$ be a non-empty class and let $\psi_A$ be a function such that $\psi_A(w)=
\zero$ if $w\in A$ and
$\psi_A(w)\diverges$ otherwise
(the \emph{pseudocharacteristic function}). Then $A$ is called \emph{semidecidable} if $\psi_A$ 
is computable. 
Fixing some 
$\Omega\subseteq D$ representing the relevant time of the model, we say that $A$ is \emph{computably enumerable} if there is a program $p\in P$
such that $A = \{I(T,p)(w)\,;\,w\in\Omega\}$.
\end{defn}

\begin{thm}[Folklore]\label{thm:ce}
For the model of ordinary Turing computability and any 
non-empty set $A\subseteq\Sigma^*$, the following hold:
\begin{enumerate}[(i)]
\item the set $A$ is semidecidable if and only if
it is the range of a partial computable function $f\colon \Sigma^*\dashrightarrow\Sigma^*$ and
\item the set $A$ is computably enumerable if and only 
if it is semidecidable.
\end{enumerate}
\end{thm}

The equivalence (i) is a classical textbook argument \cite[Theorem V]{rogers1987theory}; in equivalence (ii), the forwards direction 
is a trivial consequence of (i) and the backwards direction uses the computable bijection between
$\Sigma^*$ and the relevant time $\omega$. 
So, adapting this proof to the case of symmetric computability will preserve the equivalence (i) and the forwards direction of (ii).

\begin{thm}\label{thm:tsce}
For the model of symmetric computability and any non-empty class $A\subseteq\infwords$, the following hold:
\begin{enumerate}[(i)]
\item the class $A$ is semidecidable if and only if
it is the range of a computable class function $f\colon \infwords\dashrightarrow\infwords$ and
\item if the class $A$ is computably enumerable, then it is semidecidable.
\end{enumerate}
\end{thm}

In 
comparison to Theorem \ref{thm:ce}, the converse of (ii) is missing in
Theorem \ref{thm:tsce}; it turns out that this 
difference 
is crucial for our quest for the
desired symmetry from \S\,\ref{sec:treble}. We write $\mathsf{SC}$ for the statement 
``the class $\infwords$ is symmetrically computably enumerable'',
call this the \emph{symmetry condition}, 
and note that it is a 
set theoretic sentence in the sense of \S\,\ref{sec:classtheory}.

\begin{prop}\label{prop:sc}
The symmetry condition $\mathsf{SC}$ is equivalent to the statement ``every symmetrically semi-decidable class is symmetrically computably enumerable''.
\end{prop}

\begin{proofsketch}
    Clearly, $\infwords$ is semi-decidable, so ``$\Leftarrow$'' is obvious. For ``$\Rightarrow$'', let 
$g:\cOrd\to\infwords$ be a computable enumeration
and $A$ be any semi-decidable class
By Theorem \ref{thm:tsce}, we have a computable surjection
$f:\infwords\to A$. Then $f\circ g$ enumerates $A$.
\end{proofsketch}

The symmetry condition
expresses that the classes $\infwords$ and $\Ord$ can be identified via the computable listing provided by $\mathsf{SC}$. Therefore, assuming $\mathsf{SC}$, time, space, programs, and in-/outputs can be considered the same type of object, and the
model of symmetric computability 
has the symmetry exhibited by ordinary Turing 
computability.

We shall now see that $\mathsf{SC}$ is independent from $\GB$ and characterise under which circumstances $\mathsf{SC}$ holds.
A crucial ingredient to prove our characterisation is the following result which is a straightforward relativisation of {\cite[Theorem 6.2]{koepke2005}}.

\begin{thm}[Koepke]
\label{thm:koepke}
Let $x$ be a set of ordinals. Then any $w\in\infwords$ is in $\mathrm{L}[x]$ if and only if there is a finite program $p$ and 
$v\in\cOrd$ such that the Koepke computation of $p$ with input $v$ and parameter $x$ halts and produces the output~$w$.
\end{thm}

\begin{proofsketch}  
The backwards direction follows from the fact that a Koepke computation from a parameter $x$ is absolutely defined. Thus, if a Koepke machine produces the output $w$ upon input $x*v$, then $w$ lies in every model containing both $x$ and $v$. Since $v\in\cOrd\subseteq\mathrm{L}[x]$, we have $w\in\mathrm{L}[x]$.
For the forwards direction, assume 
$w\in\mathrm{L}[x]$ and let $\alpha$ be an exponentially closed ordinal such that $w\in\mathrm{L}_\alpha[x]$. Then by 
\cite[Theorem 7 (a)]{KoepkeSeyfferth2009}, 
$w$ is $\alpha$-Koepke computable, and thus 
Koepke computable from the parameter giving $\alpha$, i.e., Koepke computable from a parameter in $\cOrd\subseteq\FS$.
\end{proofsketch}

\begin{lem}\label{lem:VeqLSO}
If $x$ is a set of ordinals and $\infwords\subseteq\mathrm{L}[x]$, then $\mathrm{V}{=}\mathrm{L}[x]$. 
\end{lem}

\begin{proofsketch}
Assume that $\infwords\subseteq \mathrm{L}[x]$. Assume by contradiction that $\mathrm{V}\neq \mathrm{L}[x]$. Let $A$ be an $\in$-minimal set not in $\mathrm{L}[x]$, i.e., 
$A\notin \mathrm{L}[x]$, but $A\subseteq \mathrm{L}[x]$. There is a bijection $G:\Ord\rightarrow \mathrm{L}[x]$ definable from $x$, (cf., e.g., \cite[p.\ 193]{jech}). Define $w(\alpha){\downarrow}=\zero$ if and only if $G(\alpha)\in A$; then $w\in\infwords\subseteq\mathrm{L}[x]$. But then $A=\{ G(\alpha)\in \mathrm{L}[x] ; w(\alpha)=\zero\}$ and therefore $A\in\mathrm{L}[x]$.
\end{proofsketch}

\begin{thm}\label{thm:SC implies V=L}
The symmetry condition $\mathsf{SC}$ 
holds if and only if 
the universe is constructible from a set of ordinals.
\end{thm}
\begin{proofsketch}
For ``(ii)$\Rightarrow$(i)'', use the (computable) Gödel pairing function to get a computable bijection $C:\cOrd\to \cOrd\times\FS$
and identify the finite programs with $\FS$. If $u\in\cOrd$, let $C(u)= (v,p)$, and let $F(u)$ be the result of running the finite program $p$ on input $v$. By Theorem \ref{thm:koepke}, $F$ enumerates $\infwords$.

For ``(i)$\Rightarrow$(ii)'', assume $\mathsf{SC}$. 
Let $p$ be the program of the computable enumeration of $\infwords$ (which can be encoded as a set of ordinals).
By Lemma \ref{lem:VeqLSO}, it is enough to show that $\infwords\subseteq \mathrm{L}[p]$. But this follows from the fact that $p$ defines a class surjection from $\Ord$ onto $\infwords$.
\end{proofsketch}

It follows from Theorems \ref{thm:GBC_V=L} \& \ref{thm:SC implies V=L} that $\mathrm{SC}$ is independent from $\GBC$. 
We note that in the special case of $\VeqL$,
Koepke computability and symmetric computability
are equivalent (cf.\ \cite[Exercise 3.9.7 (d)]{carl}); however, if we take any nonconstructible set of ordinals $z$, then 
$\mathrm{L}[z]$ is a model of $\mathsf{SC}$ by
Theorem \ref{thm:SC implies V=L}, the set $z$ is 
symmetrically computable (by Proposition \ref{prop:setsarecomputable}), but not Koepke computable by Theorem \ref{thm:koepke} (letting $x=\varnothing$), so the two models of computability are different. This also answers \cite[Question 5.12]{lewis} about separating the stronger versions $\mathsf{SC}_\kappa$ (``$\infwords$ is computably enumerable by a program of size ${<}\kappa$''):
e.g., if $x$ is a non-constructible real, then 
$\mathsf{SC}_{\aleph_1}$ holds in $\mathrm{L}[x]$, but not $\mathsf{SC}_{\aleph_0}$.

\bibliographystyle{abbrv}
\bibliography{cie2022}
\end{document}